\documentclass[12pt,a4paper]{article}
\usepackage{graphicx,amssymb,amsmath,amsthm}
\usepackage[noadjust]{cite}
\usepackage{graphicx}
\newtheorem{theorem}{Theorem}
\newtheorem{lemma}{Lemma}

\newcommand{\be}{\begin{equation}}
\newcommand{\ee}{\end{equation}}
\newcommand{\ben}{\begin{eqnarray}}
\newcommand{\een}{\end{eqnarray}}
\newcommand{\ra}{\rangle}
\newcommand{\la}{\langle}

\newcommand{\til}{\tilde}

\newcommand{\h}{{\cal{H}}}

\newcommand{\sumn}{\sum_{n=1}^N}
\newcommand{\sumnj}{\sum_{\substack{n=1\\  n \ne j}}^N}

\newcommand{\sep}{\,;\ \ }
\newcommand{\spa}{\, ; \,}
\newcommand{\Nj}{N-1}
\newcommand{\rj}{R_j}
\newcommand{\vn}{V_{N}}
\newcommand{\vnj}{\til{V}_{N-1}}

\newcommand{\pvn}{\hat{P}_{V_N}}

\newcommand{\fvnj}{f_{\vnj}}

\newcommand{\altt}{{\beta}}
\newcommand{\kaln}{\alpha_n}
\newcommand{\aln}{\alpha_n}

\newcommand{\alj}{\alpha_j}
\newcommand{\kalj}{\alpha_j}
\newcommand{\kaltnN}{\altt_n^N}
\newcommand{\kaltjN}{\altt_j^N}
\newcommand{\kaltnNj}{\altt_n^{\Nj}}

\newcommand{\cnN}{c_n^{N}}
\newcommand{\cjN}{c_j^{N}}
\newcommand{\cnNj}{c_n^{\Nj}}

\newcommand{\subnj}{ \;;\;n=1,\ldots, j-1, j+1, \ldots,N}
\newcommand{\subnjw}{n=1,\ldots, j-1, j+1, \ldots,N}
\newcommand{\subn}{\;;\;n=1,\ldots,N}
\newcommand{\spann}{{\mbox{span}}}

\newcommand{\faccj}{\frac{\la\kaltnN,\kaltjN\ra}{||\kaltjN||^2}}

\newcommand{\alkk}{\alpha_{{k+1}}}

\newcommand{\belnkk}{\beta_{{n}}^{k+1}}
\newcommand{\belnk}{\beta_{{n}}^{k}}

\newcommand{\psikktt}{\frac{{\psi}_{k+1}}{||{\psi}_{k+1}||^2}}

\begin{document}
\baselineskip = 1.7\baselineskip

\title{Backward Optimized Orthogonal Matching Pursuit 
Approach}
\author{Miroslav Andrle, Laura Rebollo-Neira$^\ast$, 
and Evangelos Sagianos\\
NCRG, Aston University,\\
http://www.ncrg.aston.ac.uk\\
Birmingham B4 7ET,\\
United Kingdom}
\date{}

\maketitle
\begin{abstract}
A recursive approach for 
shrinking coefficients of an atomic decomposition
is proposed. The corresponding algorithm evolves 
so as to provide at each iteration a)
the orthogonal projection of a signal onto 
a reduced subspace and b) the index of the coefficient to be 
disregarded in order to construct a coarser 
approximation minimizing the norm of the 
residual error.
\end{abstract}
\noindent
EDICS Category: 1-TFSR.
\section{Introduction}
Adopting the terminology early introduced in \cite{mz}
we call atomic decomposition of a signal to the 
linear expansion:
\be
f= \sumn c_n \aln,
\label{atde}
\ee
where the atoms $\aln$ are elements of a non-orthogonal 
sequence in  the space of the signal $f$, which 
is assumed to be a Hilbert space $\cal H$. Within 
the general Matching Pursuit (MP) framework 
\cite{mz,pati,dav,do,relo} 
the atoms are chosen, by different criteria, 
from a in general redundant set which is called a 
dictionary. 
The problem of selecting  atoms in order to 
construct the signal representation with the 
minimum possible number of coefficients is a very 
complex problem. In particular, the problem of 
constructing the optimal approximation with $N$-atoms 
selected so as to  
minimize the Euclidean distance between 
the signal and the approximation is a NP-hard problem 
\cite{dav}. Thus, in this line only suboptimal 
solutions are actually feasible. 
In a previous publication a suboptimal iterative 
pursuit strategy, 
which is only optimal
at each iteration step, has been introduced 
with the name of Optimized Orthogonal Matching
Pursuit (OOMP) approach \cite{relo}. Such a technique  
is based on an iterative method for computing
dual  atoms $\kaltnN \subn$ giving rise
to a representation
of the orthogonal projector operator onto the
subspace generated
by the set of atoms $\kaln \subn$. These
atoms, which happen to be biorthogonal to
atoms $\kaln \subn$ \cite{biofo},  
allow to compute the coefficients in
(\ref{atde}) by computing inner products, i.e.
$c_n^N= \la \kaltnN , f \ra$
(the superscript $N$ indicates that the
dual atoms, and therefore the coefficients,
depend of the number $N$ of dictionary atoms being
considered). The  OOMP approach tackles the problem 
of selecting the 
new atom $\alpha_{N+1}$ to improve the approximation. 
Moreover, the  coefficients of
the atomic decomposition are  recursively
modified in order to yield an 
optimal approximation in the enlarged subspace. 
Here we consider the reverse situation: 
We assume that an atomic decomposition is given and
we wish to eliminate some coefficients.
For such an end we propose a technique, that we term 
Backward Optimized Orthogonal Matching Pursuit (BOOMP)
which consists of the following elements: 
a) a recursive approach to modify the coefficients of 
the atomic decomposition when one of the coefficients 
is to be disregarded and b) the criterion to select such  
a coefficient. 
Although the technique can be applied to reduce
coefficients of any atomic decomposition, regardless
of how such a decomposition is obtained, in this
letter we focus on its implementation
as a posterior step of OOMP. The
reason for taking this route is that the
implementation of BOOMP is really straightforward
on the outputs of OOMP. In addition, we believe that
this way of introducing the approach makes
more clear its purpose and also
its implementation steps.

We would like to stress that even the 
construction of suitable suboptimal signal approximations 
by using dictionaries is a 
complex theoretical and practical 
challenge \cite{devo,ti}. 
This communication aims at 
enhancing the fact that, since suboptimal forward and 
backward approximations are in general not reversible,  
application of the proposed backward pursuit approach, 
after a forward pursuit selection of atoms, may result in 
a  gain with respect to sparseness of the representation. 
This is clearly illustrated here by a simple example. 
\section{Adaptive MP strategies}
The MP approach  is a technique 
to compute adaptive signal representations by 
iterative selection of atoms \cite{mz}. In 
its original form this technique 
does not yield at each iteration the linear
expansion of the selected atoms that approximate the
signal at best in a minimum distance sense.
A later refinement, which does provide such an
approximation, has been termed Orthogonal Matching Pursuit 
(OMP) \cite{pati,dav,ma}. However, since
 OMP selects the atoms according to the MP prescription,
the selection criterion is not optimal in 
the sense of minimising the residual of the 
new approximation. The requirement  of such minimization 
has led to the recently introduced 
OOMP approach \cite{relo}. 
This technique is implemented by means of an adaptive  
biorthogonalization method  which, within the 
workings of the selection process, generates  the
set of biorthogonal atoms yielding 
orthogonal projections \cite{relo,biofo}. 
Such biorthogonal atoms 
are used to compute the coefficients of the 
atomic decomposition and 
are obtained through the following recursive equations: 
\begin{align}
\belnkk&= \belnk - 
\beta_{{k+1}}^{k+1} \la \alkk, \belnk  \ra
\sep n=1,\ldots,k,
\label{e1}\\
\beta_{{k+1}}^{k+1}& = \psikktt, 
\label{e2}
\end{align}
where the set $\psi_j\spa j=1,\ldots,k+1$  
is obtained inductively, from  $\psi_1=\alpha_{1}$, 
by orthogonalization 
of atoms $\alpha_{j}\spa j=1,\ldots,k+1$.  
These atoms are selected from 
the dictionary by minimizing, at each iteration step, 
the norm of the residual error in approximating the 
signal \cite{relo}.

Let us stress once more that, since all the above 
mentioned pursuit strategies evolve by fixing the 
atoms selected in the previous steps, there is 
plenty of room for possible improvement with regard to 
compression of the representation.

Improving compression after the 
OOMP procedure implies having to eliminate some
coefficients of the atomic resulting decomposition. 
For the coarser approximation to be optimal 
in a minimum distance sense, the remaining 
coefficients must be recalculated \cite{ze1,ze2,bioba}.
This feature of non-orthogonal expansions is a 
major difference with orthogonal ones 
and has motivated an adaptive approach 
to modify biorthogonal atoms in order for 
then to yield orthogonal projections 
when the corresponding subspace is reduced \cite{bioba}.

Let us suppose that OOMP has selected $N$ atoms to 
represent a given signal 
up to some predetermined precision.  Let us 
denote $V_N$ to the subspace spanned by such 
atoms i.e,  $\vn=\spann\{\alpha_{1}, \ldots, 
\alpha_{N}\}$ and let 
$\vnj$ be the subspace
 which is left by removing one atom, 
say the $j$-th one, 
i.e. $\vnj=\spann\{\alpha_{1},\ldots,\alpha_{{j-1}},
 \alpha_{j+1}, \ldots, \alpha_N\}$.
Since the  biorthogonal atoms 
 $\beta_n^N \spa n=1,\ldots,N$ are available
(as output of the 
OOMP procedure), 
to construct the orthogonal projector
of $f$ onto $\vnj$ we just need to 
modify the  atoms $\beta_n^N \spa n=1,\ldots,j-1,j+1,\ldots,N$  
as follows \cite{bioba}
\be
\kaltnNj
= \kaltnN -
\frac{\kaltjN \la \kaltjN,\kaltnN\ra }{|| \kaltjN||^2}
\sep \subnjw. 
\label{reba}
\ee
In writing the above equation we have re-defined
the superscript $N-1$. Now this
upper index indicates that the biorthogonal atoms
are modified in order to account for
the deleting of any one atom (not necessarily the  
last element of the  spanning set).
In the next section we discussed how these adaptive
backward equations generate the proposed
BOOMP approach.
\section{Backward Optimized Orthogonal Matching Pursuit} 
After the selection of $N$ atoms the 
OOMP approach provides a representation of 
a signal $f$ as given by \cite{relo}
\be
f_{\vn}= \pvn f= \sumn c_n^N \kaln,
\label{five}
\ee
$\pvn f$ indicates the orthogonal projection of 
the signal $f$ onto $V_N$ and 
coefficients $c_n^N$ are  
obtained as
$c_n^N= \la \kaltnN ,f\ra$. 
Theorem 1 below proves that, if we 
decide to eliminate the coefficient 
$c_j^N$ from the above expansion, in 
order to obtain the optimal approximation of  
$f$ in the reduced subspace $\vnj$, 
the remaining coefficients 
$c_n^N \spa n= 1,\ldots,j-1,j+1,\ldots,N$ 
should be modified as follows: 
\be
\cnNj= \cnN - \faccj \cjN \sep n=1,\ldots,j-1,j+1,\ldots,N. 
\label{cnNj}
\ee 
For the sake of organizing the corresponding proof let us 
first prove the following lemma:
\begin{lemma}\label{lem1}
 Let signal $\fvnj$ be given by 
\be
\fvnj = \sumnj \cnNj \kaln, 
\label{fvnj}
\ee
with coefficients 
$\cnNj$ as in (\ref{cnNj}), and let $f$ be a signal in $\h$.
The difference $f-f_{\vnj}$ is orthogonal to every function 
in ${\vnj}$.
\end{lemma}
\begin{proof} 
Using (\ref{cnNj}) and (\ref{fvnj}) we have: 
\begin{align}
f-f_{\vnj}&=f - \sumnj \cnNj \kaln = 
                f - \sumnj \cnN \kaln + 
               \sumnj \faccj \cjN \kaln \nonumber\\
 &= f - \sumn \cnN \kaln + \cjN \kalj  + \sumn   \faccj \cjN \kaln - \cjN  \kalj. 
\label{dif}
\end{align}
Since $\kaltjN \in \vn$, it follows from 
(\ref{five}) that 
$\sumn \faccj \kaln =\frac{ \pvn\kaltjN}{||\kaltjN||^2} = \frac{\kaltjN}{||\kaltjN||^2}$. Hence (\ref{dif}) 
turns out to be
\be\label{dif2}
f-f_{\vnj} =  f - \sumn \cnN \kaln + \cjN  \frac{\kaltjN}{||\kaltjN||^2}.
\ee
Since by hypothesis $\sumn \cnN \kaln = \pvn f$, 
the difference $f - \sumn \cnN \kaln$ is orthogonal to every 
function in $\vn$. Furthermore, since  
$\la \kaln , \kaltjN \ra =\delta_{n,j}$,  by taking 
the inner product both sides of (\ref{dif2}) with every function 
$\kaln\subnj$ we obtain $\la \kaln ,f-f_{\vnj} \ra = 0$, 
which proves that $f-f_{\vnj}$ is orthogonal to every function 
in $\vnj$.
\end{proof}
\begin{lemma}\label{orthproj} 
The coefficients $\cnNj$  of the linear expansion 
\be
\fvnj =  
\sumnj \cnNj \kaln 
\ee
minimizing the distance in $\vnj$ to a given signal 
$f \in \h$ are obtainable, from  
 $\cnN$ and    
 $\kaltnN \subn $, as prescribed in
(\ref{cnNj}).
\end{lemma}
\begin{proof}
 Let $g$ be an arbitrary signal in $\vnj$ and let us
write $|| f - g||^2$ as follows: 
\be
|| f - g||^2 = || f - f_{\vnj}  + f_{\vnj} -g ||^2.
\ee
From Lemma~\ref{lem1} we know that $f - f_{\vnj}$ is orthogonal 
to every function in $\vnj $ and since 
 $f_{\vnj}-g$ is  
in $\vnj$ we have
\be
||f - g||^2 = ||f - f_{\vnj}||^2 + ||f_{\vnj} -g||^2, 
\ee
from where we conclude that $||f - g||^2$ is minimized if 
$g \equiv f_{\vnj}$.
\end{proof}
Lemma~\ref{orthproj} tells us how to proceed to 
disregard coefficients of a non-orthogonal linear expansion. 
Assuming that the 
coefficient $c_j^N$ to be disregarded has been selected, 
in order to optimize the approximation in a minimum distance 
sense, the remaining coefficients should be modified as 
indicated in (\ref{cnNj}).
The next theorem  gives an answer to the question as to how 
to select the coefficient $c_j^N$ to be neglected.
\begin{theorem}\label{maintheo}
 Let $\rj$ be the residual resulting 
by disregarding a coefficient $c_j^N$ for passing from 
approximation $f_{\vn}$ to $f_{\vnj}$ i.e., 
$f_{\vn}= f_{\vnj} + \rj$. 
In order to minimize the norm of the residual $\rj$ such 
coefficient is to be chosen as the one yielding a minimum 
value of the quantity
\be
\frac{|\cjN|^2}{||\kaltjN||^2}.
\label{cond}
\ee
\end{theorem}
\begin{proof}
 Since $\rj =f_{\vn} - f_{\vnj} = 
\sumn \cnN \kaln - \sumnj \cnNj \kaln $, by using 
(\ref{cnNj}) we have: 
\begin{align}
\rj&=\sumn \cnN \aln -\sumnj \cnN \aln + \sumnj 
\cjN \aln \faccj\nonumber\\
 &=\cjN \kalj + \cjN \sumn \kaln \faccj - \cjN \alj
\end{align}
As already discussed, 
$\sumn \aln {\la \kaltnN , \kaltjN \ra}= \kaltjN$.  Then,    
from the last equation it follows that 
$\rj= \frac{\cjN \kaltjN}{||\kaltjN||^2}$. 
Consequently, 
in order to 
minimize $||\rj||^2$ the coefficient $\cjN$ to be 
neglected is the one minimizing (\ref{cond}). 
\end{proof}
Theorem~\ref{maintheo} leads to a recursive algorithm for shrinking  
coefficients.
We call such algorithm BOOMP, because, 
at each iteration, it selects the atom to be deleted  
according to a selection criterion which is equivalent to 
the one proposed by OOMP \cite{relo} for
forward approximations. Moreover, BOOMP is a natural 
complement of OOMP because its implementation on 
the output of OOMP is extremely simple. The few 
necessary steps are describe below.
{\subsection*{BOOMP algorithm}}
Let us assume that atoms $\aln \subn $  have been selected by 
the OOMP approach in order to approximate a signal $f$. 
Hence, the biorthogonal set $ \kaltnN $ and the 
corresponding coefficients $c_n^N \subn$ are also known. 
The BOOMP approach for reducing coefficients evolves as follows: 
\begin{itemize}
\item
Select the index $j$ of the coefficient $\cjN$ 
to be disregarded as the one yielding a minimum value of 
the quantity $ \frac{|\cjN|^2}{||\kaltjN||^2}$ as 
$j$ ranges from $1$ to $N$.
\item
Modify the corresponding biorthogonal atoms and coefficients 
as prescribed in (\ref{reba}) and (\ref{cnNj}) respectively.
\item
 Set $N=N-1$ and repeat the above steps until the coarsest acceptable 
approximation is reached.
\end{itemize}
{\subsection*{Example}
We illustrate now  by a simple example
the main remark
of this communication: namely that BOOMP can
improve the compression performance of the forward
OOMP approach.
We construct a dictionary of Mexican hat wavelets
 given by the functions
\be
\alpha_{m,n}(t)= 2^{\frac{m}{2}}\alpha(t 2^{m} - 0.2n)\\ 
\;\;\;\;\;\;\text{ with } \;\;\;\;\;
\alpha (t)=\frac{2}{\sqrt{3}}\pi^{-\frac{1}{4}}
 (1-t^{2}) e^{- \frac{t^{2}}{2}}.
\ee
By considering  scales $m=0,1,2,3$, and 4  to
cover the $[0,4]$ interval we have  a dictionary of 
$665$ atoms.
The signal $f$ to be represented is a chirp generated 
by the MATLAB instructions: 
\begin{verbatim}
t=0:0.01:4; f=chirp(t,0,1,2); 
\end{verbatim}
In order to produce a high quality representation of this  
chirp the OOMP approach
selects $N=60$ atoms (the norm of the residual error 
is 0.0544). By applying the proposed
BOOMP to the OOMP approximation we reduce the
number of coefficients up to $34$ and the approximation
is the one depicts in the top left graph of Figure~1.
The norm of the residual error with respect to the 
true signal is in this case 1.18. 
However, if rather than applying the BOOMP approach
we stop the OOMP approach at iteration $34$,  so as
to have the same number of coefficients as in
the previous case, the approximation is the
one shown in the  top right graph of Figure~1 
and the  norm of the residual error is $1.77$. It is clear from 
the graphs that,  in addition to yielding the smallest  
residual error with the same number of coefficients, 
the approximation  obtained by the 
BOOMP approach is overall visually superior to the  
OOMP one with the same number of coefficients. 
The bottom left and right graph depict, respectively, 
the absolute value of the difference between the 
chirp signal and the corresponding approximations. 
\section{Conclusions}
A recursive approach for shrinking coefficients 
of an atomic decomposition
has been proposed. The approach is based on an 
adaptive technique which 
allows to modify biorthogonal functions in order to 
yield orthogonal projectors onto a reduced subspace. 
A criterion for disregarding coefficients has been 
discussed. Such criterion leads to an 
iterative procedure that we have termed BOOMP, because 
it evolves so as to fulfil identical 
requirements to those of the OOMP method. 
Accordingly, BOOMP  provides 
at each iteration a) the coefficient 
of the atomic decomposition to be  deleted  
in order to construct a coarser approximation 
minimizing the norm of the residual error 
b) the coefficients of such decomposition   
rendering optimal approximation in the 
same sense. The approach is a good complement 
to  OOMP, very simple to implement and definitely 
worth trying in any case. However it is appropriate 
to stress  that situations for which the 
forward approach renders better results that 
the  combination with the backward one 
should certainly exist. This is a consequence of 
lack of global optimality in both directions.
Finally we would like to remark that if, 
rather than (\ref{cond}), one decided to
apply another criterion for
disregarding coefficients
(see \cite{deno2,deno3} for some alternative ones) 
in order to leave an approximation minimizing the 
 distance to the signal the 
remaining coefficients should be modified as
prescribed in (\ref{cnNj}). 

\section*{Acknowledgements}
We would like to thank two anonymous referees for their 
valuable comments on a previous submittal.

MATLAB codes for implementation of both OOMP and 
BOOMP are available upon request.
 
Support from EPSRC (GR$/$R86355$/$01) is acknowledged.

\newpage
\begin{figure}[ht]
\centering
\includegraphics[angle=0,width=7cm]{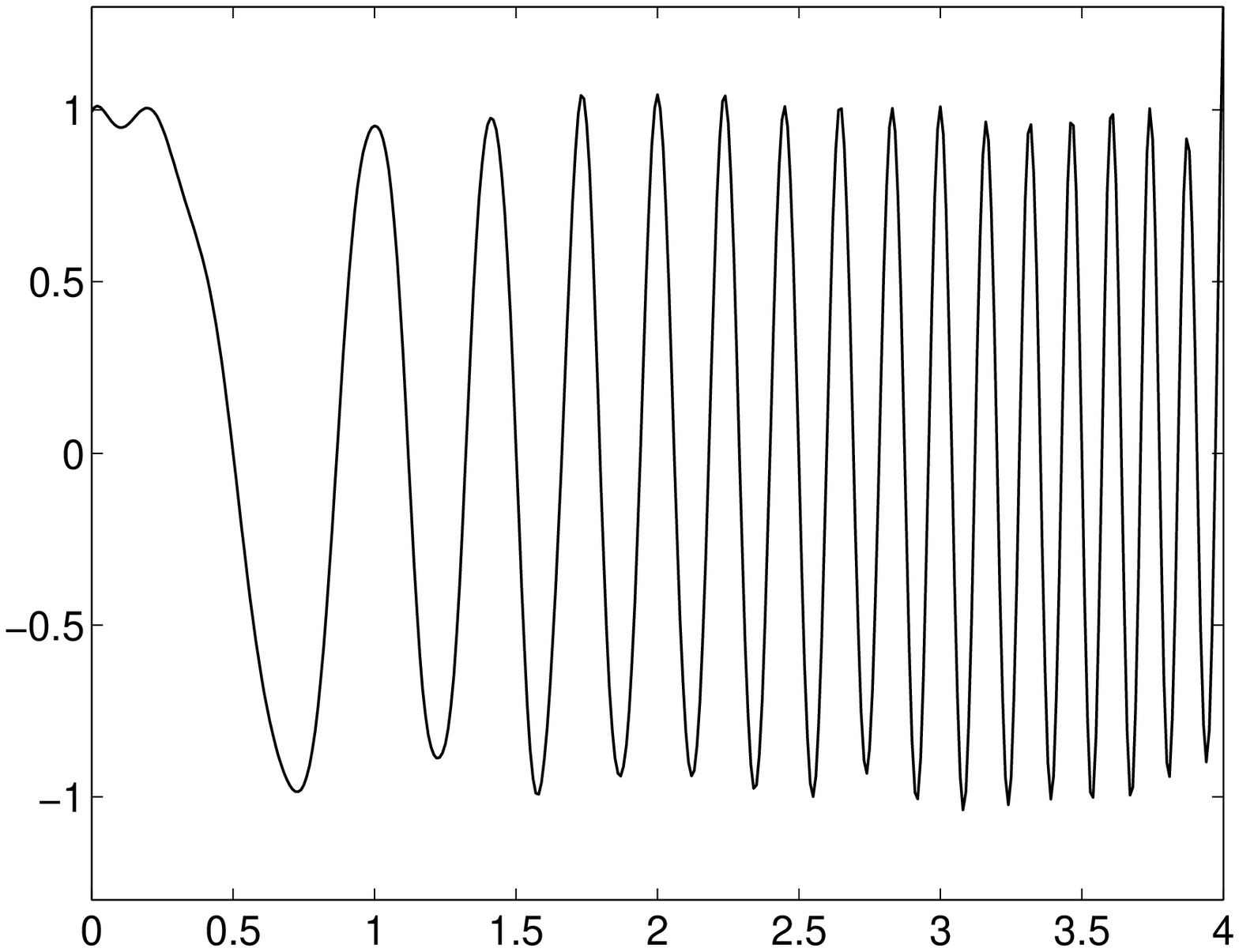}\qquad
\includegraphics[angle=0,width=7cm]{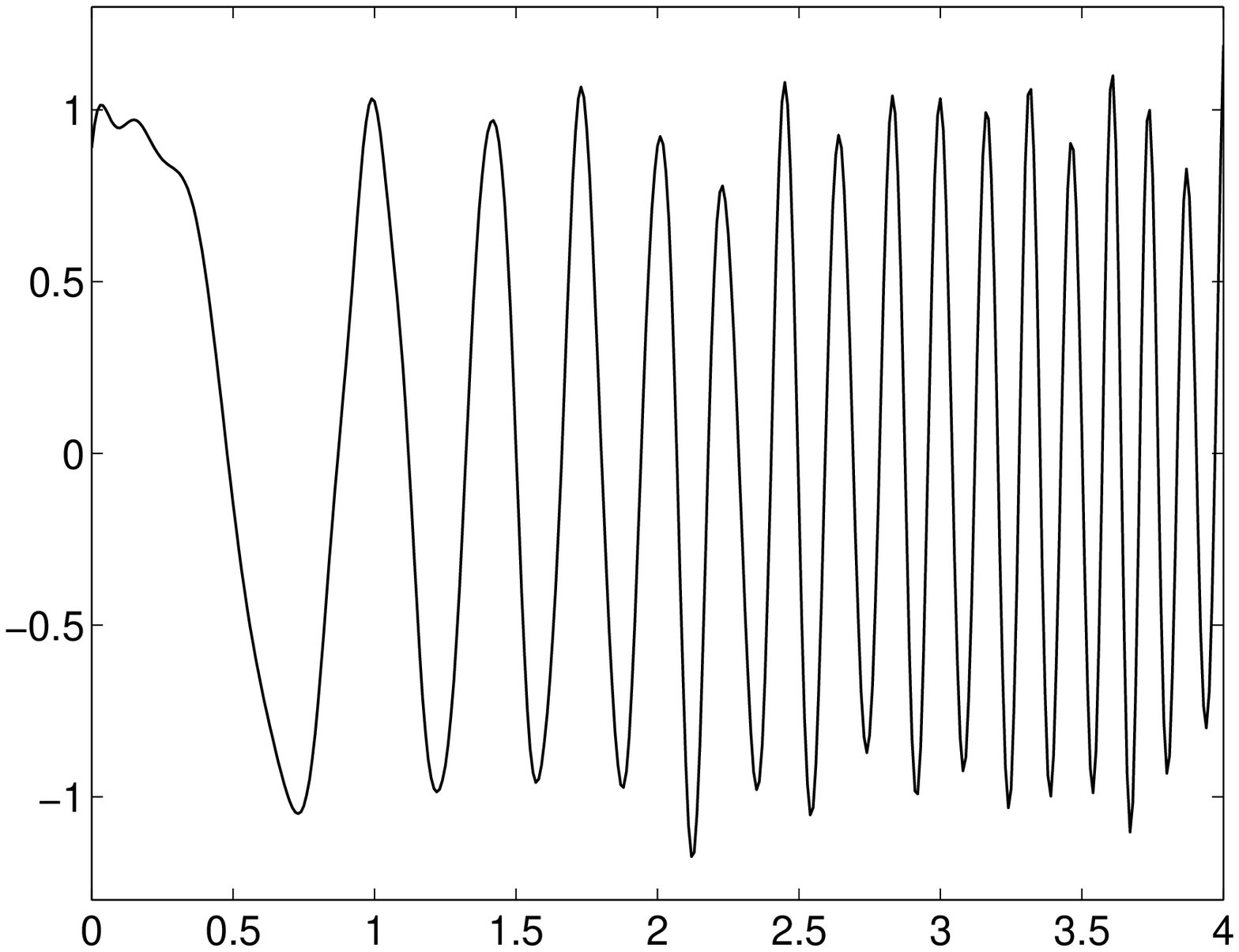}\\[0.5cm]
\includegraphics[angle=0,width=7cm]{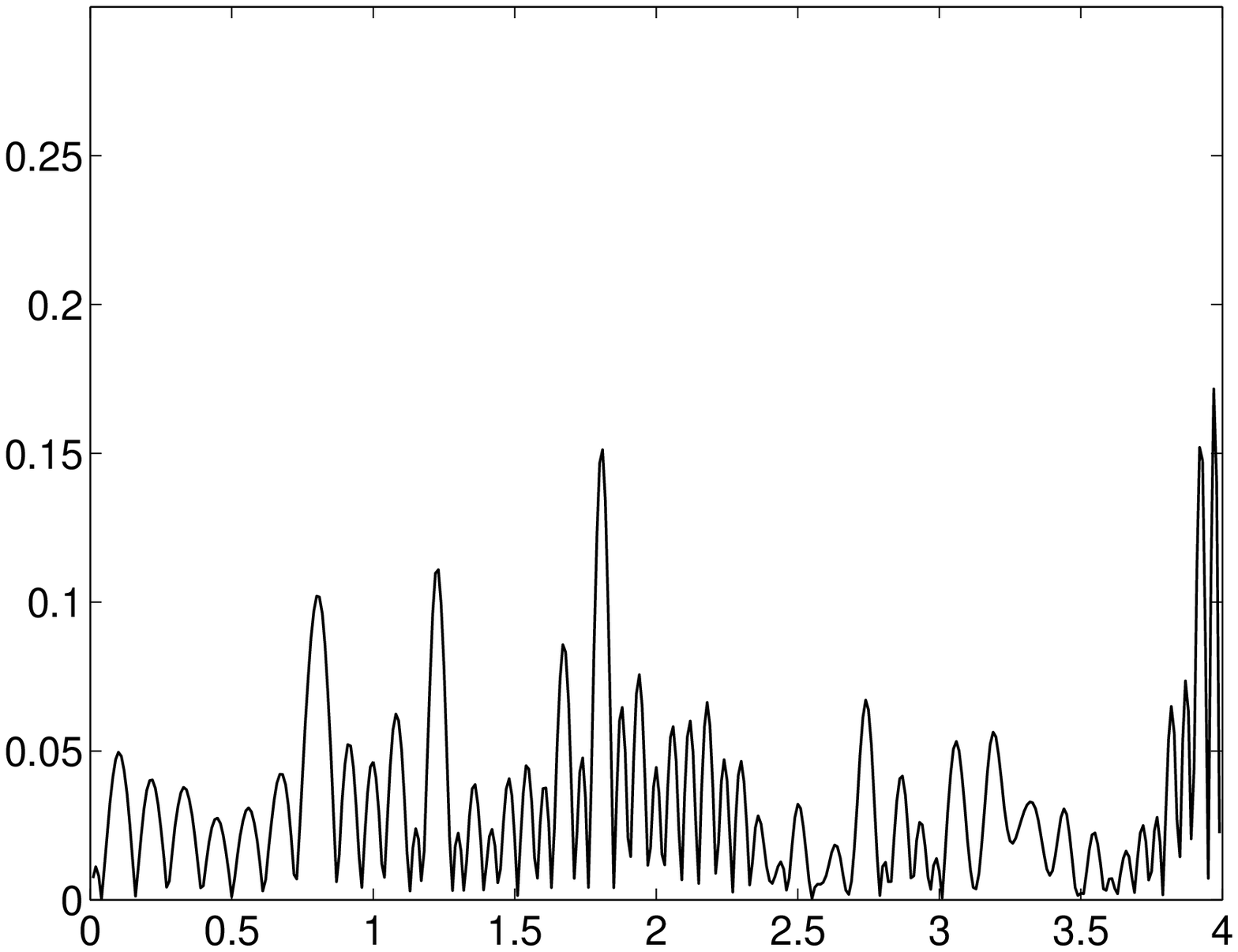}\qquad
\includegraphics[angle=0,width=7cm]{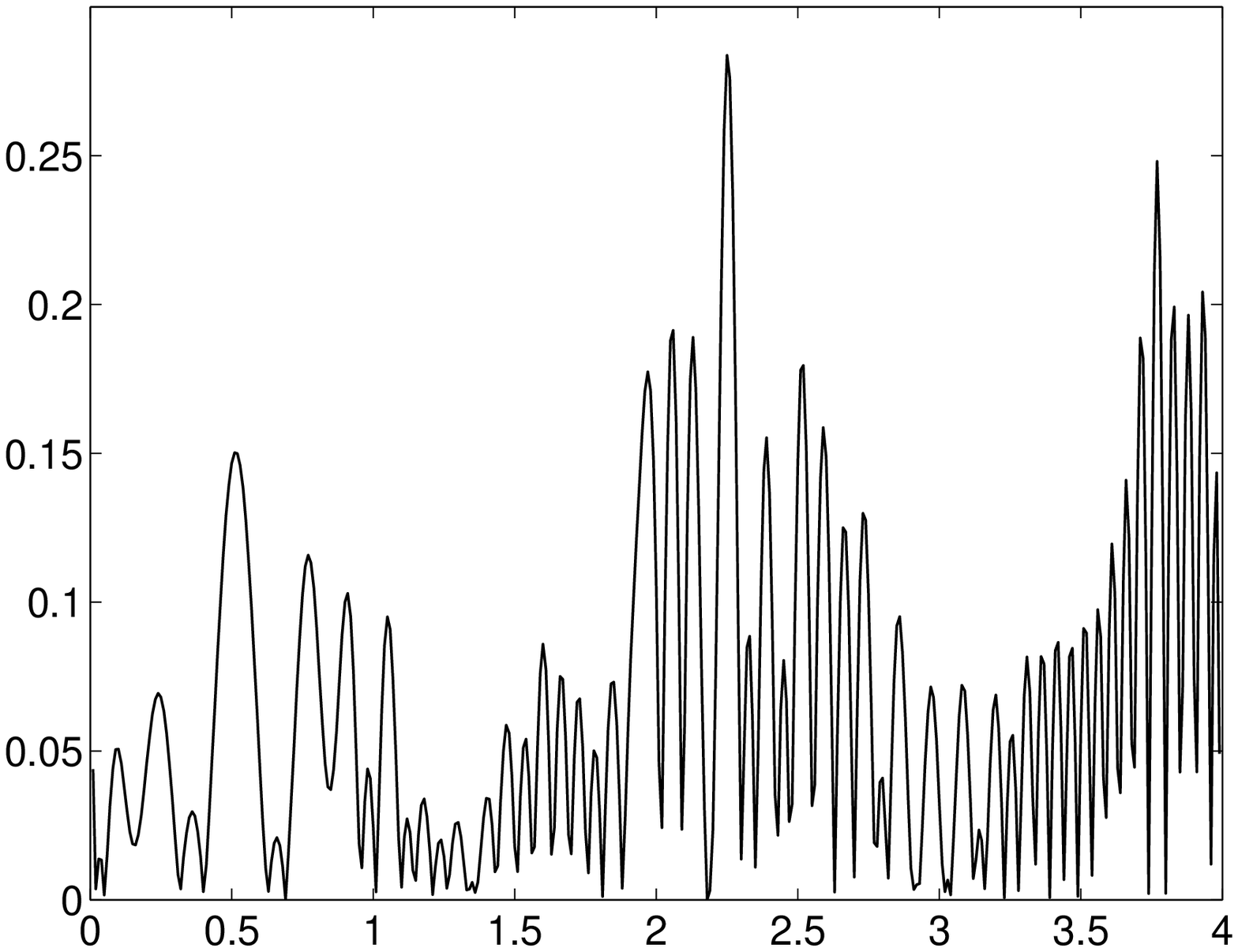}
\label{differ}
\caption{The top left graph represents the approximation
obtained by applying the BOOMP approach for
reducing the 60 coefficients of the high quality
OOMP approximation up to 34. The top right
graph corresponds to the OOMP approximation
resulting by stopping the approach after the
selection of 34 atoms. The bottom left
graph depicts the absolute value of the difference  between
the chirp signal and the BOOMP approximation
of the left top graph. The right
bottom graph has the same description as the
left one, but with respect to the approximation
of the right top graph.}
\end{figure}
\newpage
\end{document}